\documentclass[12pt]{article}

\usepackage{amsmath}
\usepackage{epsfig}                     
\usepackage{amssymb}
\usepackage{fullpage}

\newtheorem{thm}{Theorem}
\newtheorem{cor}[thm]{Corollary}
\newtheorem{defn}{Definition}

\newtheorem{lemma}[thm]{Lemma}

\newtheorem{prop}[thm]{Proposition}

\newenvironment{proof} { \emph{Proof.} } { {\rule{2mm}{2mm}}\\}

\newcommand{\bea}{\begin{eqnarray*}}
\newcommand{\eea}{\end{eqnarray*}}

\newcommand{\e}{\epsilon}

\newcommand{\R}{\mathbb{R}}

\newcommand{\Z}{\mathbb{Z}}
\newcommand{\Q}{\mathbb{Q}}

\renewcommand{\phi}{{\varphi}}

\newcommand{\bm}{\begin{pmatrix}}
\newcommand{\fm}{\end{pmatrix}}

\begin{document}
\title{Bilipschitz equivalence is not equivalent to quasi-isometric equivalence for finitely generated groups}
\author{Tullia Dymarz}
\maketitle
\begin{abstract} 
We show that certain lamplighter groups that are quasi-isometric to each other are not bilipschitz equivalent. This gives a positive answer to a question in \cite{H} page 107.
\end{abstract}
A \emph{quasi-isometry} between metric spaces is a map $f: X \to Y$ such that for some $K,C>0$
$$ -C + \frac{1}{K} d(x,y) \leq d(f(x),f(y)) \leq K d(x,y) +C $$
for all $x,y \in X$ and such that $nbhd_{C}(f(X))=Y$.
This is a generalization of the more common notion of a \emph{bilipschitz equivalence}: a bijection between metric spaces that satisfies for some $K$
$$\frac{1}{K} d(x,y) \leq d(f(x),f(y)) \leq K d(x,y).$$
A natural question to ask is for which classes of metric spaces are these two notions equivalent. In \cite{Kl}, Burago and Kleiner gave an example of a separated net in $\R^2$ that is not bilipschitz equivalent to the integer lattice (but the two are quasi-isometric). Our interest is in the class of finitely generated groups equipped with word metrics. For a finitely generated group $\Gamma$ a choice of generating set $S$ determines a Cayley graph $\Gamma_S$ with metric $d_S$. The metric $d_S$ depends on $S$ but for any given group all Cayley graphs are bilipschitz equivalent.
The example in \cite{Kl} is not a Cayley graph of a finitely generated group.
In this paper we prove the following Theorem:
\begin{thm}\label{thebigone}
Let $F$ and $G$ be finite groups with $|F|=n$ and $|G|=n^k$ where $k>1$ and $k\in\Z$. Then there does not exist a bijective quasi-isometry between the lamplighter groups $G \wr \Z$ and $F \wr \Z$  if $k$ is not a product of prime factors appearing in $n$. 
\end{thm}

For discrete spaces, a bijective quasi-isometry is the same as a bilipschitz map. 
In \cite{W}, Whyte proved that for \emph{nonamenable} finitely generated groups any quasi-isometry is a bounded distance from a bijective quasi-isometry so our examples are necessarily \emph{amenable}.
Recall that a \emph{F\o lner} sequence in a discrete metric space is a sequence of finite sets $S_i$ such that for all $r>0$
$$ \lim_{i \to \infty} \frac{|\partial_r S_i|}{|S_i|} \to 0$$ 
where $\partial_r S_i$ denotes all points that are either not in $S_i$ but at a distance of at most $r$ from $S_i$ or points in $S_i$ at distance at most $r$ from points not in $S_i$. 
A group is said to be \emph{nonamenable} if it does not admit a F\o lner sequence and \emph{amenable} otherwise.
Using uniformly finite homology, one can check when a map between (uniformly discrete bounded geometry) spaces is bounded distance from a bijection (see \cite{W} or Section \ref{UFH} for details). 
The obstruction for a quasi-isometry to be a bounded distance from a bijection can only come from the existence of a  F\o lner sequence. Since nonamenable groups do not have F\o lner sequences any quasi-isometry is a bounded distance from a bijection.

With respect respect to a certain generating set, the Cayley graph of the lamplighter group $F \wr \Z$ where $|F|=n$ is the Diestel Leader graph $DL(n,n)$.
We will describe the geometry of Diestel Leader graphs and construct explicit F\o lner sequences in these graphs in Section 
\ref{geomDL}. 
The main resource we have for analyzing quasi-isometries of Diestel Leader graphs is the following theorem of Eskin-Fisher-Whyte:
\begin{thm}\label{efw}\cite{EFW} 
Any quasi-isometry $\phi: DL(n,n) \to DL(n,n)$ is a bounded distance from a ``standard"  map of the form $(x,y,t) \mapsto (\phi_l(x), \phi_u(y), t)$ where $x,y \in \Q_n$ and $t\in \R$ and $\phi_l,\phi_u$ are bilipschitz.
\end{thm}
For an explanation of the coordinatization of $DL(n,n)$ as $\Q_n \times \Q_n \times \R$ see Section \ref{coords}.

Now since any two finite groups $G$ and $G'$ that have the same order give rise to the same Diestel Leader graph we can restrict out attention to $G=F^k$ where $F^k$ is the direct product of $k$ copies of $F$. The group $F^k \wr \Z$ appears as a finite index subgroup of  $F\wr \Z$ of index $k$ (see Section \ref{upmap}).  By \cite{D} we know that finite index subgroup inclusion $i:F^k \wr \Z \to F \wr \Z$ is not a bounded distance from a bijection. However, as pointed out in \cite{D}, if we were able to find a quasi-isometry 
$$\phi: F^k \wr \Z \to F^k \wr \Z$$ for which $|\phi^{-1}(p)|=k$, also known as a $k$ to $1$ quasi-isometry, then $i \circ \phi: F^k \wr \Z \to F \wr \Z$ would be a bounded distance from a bijective quasi-isometry.  We claim first that no such $\phi$ exists if $k$ is not a product of prime powers appearing in $n$ and second that any bijective quasi-isometry $\phi'$ between $F^k \wr \Z$ and $F \wr \Z$ arises in this way, i.e.  $\phi'$ is a bounded distance from $i \circ \phi$ for some $\phi$.

The second claim is a simple consequence of the geometry of the graph $DL(n,n)$  (see Section \ref{upmap}). The proof of the first claim proceeds by contradiction. By Theorem \ref{efw} the $k$ to $1$ map $\phi$ has the form $(\phi_l, \phi_u,id)$. Combining a theorem of Cooper \cite{FM1} on the structure of bilipschitz maps of $\Q_n$ with a key observation of Juan Souto
we are able to replace $\phi$ with a map $\bar{\phi}$ that is still $k\ to\ 1$ but that is now a bounded distance from a map $(\bar{\phi_l}, \bar{\phi_u}, Id)$ where 
$\bar{\phi_l}, \bar{\phi_u}$ are now \emph{measure linear}, i.e. on the level of measure $\phi_l, \phi_u$ scale sets by fixed constants $\lambda_l, \lambda_u$. This map $\bar{\phi}$ is constructed in Section \ref{tukiasection}.


In Proposition \ref{powersofn} we are able to argue that $\lambda_l, \lambda_u$ are products of powers of primes appearing in $n$. 
Finally in Section \ref{constructQI} we construct an explicit quasi-isometry $\psi$ that is bounded distance from $\bar{\phi}$ but that on a sequence of F\o lner sets is approximately $1/\lambda_l\lambda_u$ to $1$. Using a Theorem of Whyte \cite{W} we conclude that this is only possible if $1/\lambda_l\lambda_u=k$.

{\bf Acknowledgments.} I would like to thank Bruce Kleiner, Enrico Le Donne, Lars Louder, Amir Mohammadi, Kevin Whyte and especially Juan Souto for useful conversations. I would also like to thank Shmuel Weinberger for introducing me to this problem. 

\section{Geometry of $DL(n,n)$}\label{geomDL}
There are now many places in the literature that describe the connection between the group $F \wr \Z$ and the graph $DL(n,n)$. See for example \cite{TW,Wo}. 
We will freely borrow from these treatments. Understanding the group $F \wr \Z$ itself is not as important for us as understanding the geometry of $DL(n,n)$ so we will focus our energy on describing this graph. 

Let $T_n$ be a regular directed $n+1$ valent infinite tree where each vertex has $n$ incoming edges and one out going edge. Fixing a basepoint $p \in T_{n+1}$ we define a ``height map'' $h: T_{n+1} \to \R$ by setting $h(p)=0$ and then mapping each coherently oriented line in $T_{n+1}$ isometrically onto $\R$ with vertices mapping onto $\Z$.  Then we define
$$DL(n,n) = \{ (p,q) \in T_{n+1} \times T_{n+1} \mid h(p) + h(q) =0 \}.$$  
On $DL(n,n)$ we can also define a height map $ht: DL(n,n) \to \R$ by setting $ht(p,q)=h(p)=-h(q)$.
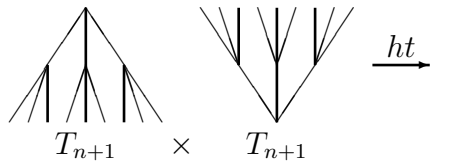
\begin{figure}[h]
\center
\setlength{\unitlength}{.1in} 
\begin{picture}(30,8)(1,.5) 
\linethickness{1pt} 

\put(5,-.25){\makebox(0,0){$T_{n+1}$}} 
\put(1,1){\line(2,3){4}}
\put(2,1){\line(1,3){1}}
\put(3,1){\line(0,1){3}}

\put(4,1){\line(1,3){1}}
\put(5,1){\line(0,6){6}}
\put(6,1){\line(-1,3){1}}

\put(7,1){\line(0,1){3}}
\put(8,1){\line(-1,3){1}}
\put(9,1){\line(-2,3){4}}

\put(10,-.25){\makebox(0,0){$\times$}} 

\put(15,-.25){\makebox(0,0){$T_{n+1}$}} 

\put(11,7){\line(2,-3){4}}
\put(12,7){\line(1,-3){1}}
\put(13,7){\line(0,-1){3}}

\put(14,7){\line(1,-3){1}}
\put(15,7){\line(0,-6){6}}
\put(16,7){\line(-1,-3){1}}

\put(17,7){\line(0,-1){3}}
\put(18,7){\line(-1,-3){1}}
\put(19,7){\line(-2,-3){4}}


\put(20,4){\vector(1,0){3}}  
\put(25,1){\line(0,1){6}}
\put(21.5,5){\makebox(0,0){$ht$}}

\end{picture}
\caption{a portion of DL(3,3)}
\label{tnxtn}
\end{figure}

\begin{figure}[h]
\begin{center}
\includegraphics[totalheight=2.2 in]{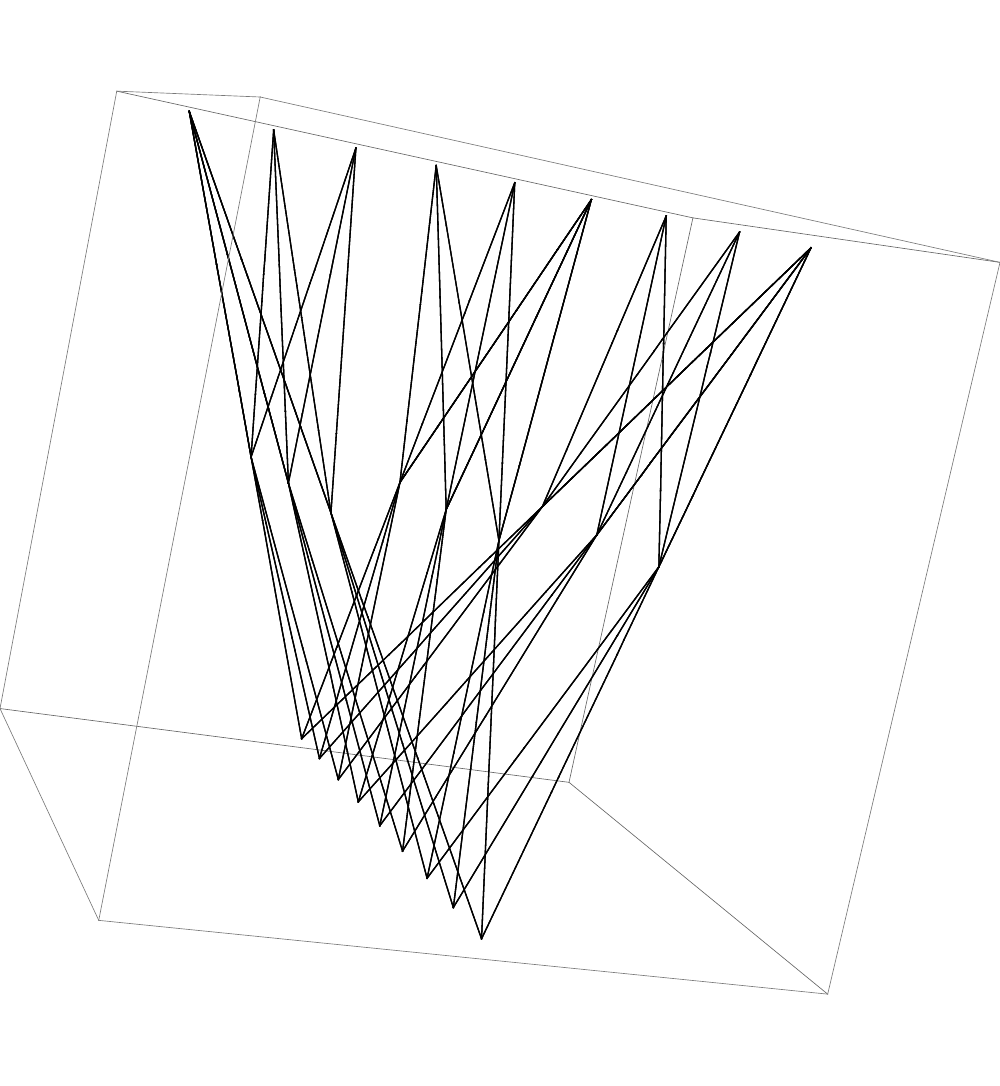}
\caption{a portion of DL(3,3)}
\label{dl33}
\end{center}
\end{figure}

By drawing the second tree ``upside down" in Figure \ref{tnxtn}, we can represent a vertex in $DL(3,3)$ by a pairs of vertices, one in each tree, that have the same height. Likewise a pair of edges at the same height represents and edge in $DL(3,3)$. A portion of the actual graph can be seen in Figure \ref{dl33}. Both figures represent a connected component of $ht^{-1}([a,b])$ where $a,b \in \Z$ with $b-a=2$.

\begin{defn}[box]\label{box} We call a connected component of $ht^{-1}([a,b])$ where $a,b \in \Z$ and $b-a=H$ a box of height $H$. 
\end{defn}
Note that  a box $S_H$ of height $H$ has volume and boundary
$$|S_H|=(H+1)n^H, \quad |\partial S_H| \sim 2 \cdot n^{H}.$$
Therefore, for any increasing sequence $H_i$, the sequence $\{S_{H_i}\}$ a F\o lner sequence.

\subsection{Subgroups of index $k$ in $F \wr \Z$}\label{upmap}

As mentioned in the introduction, the group $F^k \wr \Z$ appears as a subgroup of index $k$ in $F\wr \Z$.
To see this, let  $\rho: F \wr \Z \to \Z$ be the homomorphism that is projection onto the second factor. Then $$\rho^{-1}( k \Z) \simeq F^k \wr \Z.$$ From a geometric point of view this inclusion allows us to see $DL(n^k, n^k)$ sitting inside of $DL(n,n)$. The map $\rho$ is simply the height map $ht$. In figure \ref{dl9933}, the egdes of $DL(3^2,3^2)$ are distorted to show how $DL(3^2,3^2)$ maps into $DL(3,3)$.

\begin{figure}[h]
\setlength{\unitlength}{.1in}
\center 
\begin{picture}(20,8)(1,.5) 
\linethickness{1pt} 

\put(-4,1){\line(2,3){4}}
\put(-3,1){\line(1,2){3}}
\put(-2,1){\line(1,3){2}}

\put(-1,1){\line(1,6){1}}
\put(0,1){\line(0,6){6}}
\put(1,1){\line(-1,6){1}}

\put(2,1){\line(-1,3){2}}
\put(3,1){\line(-1,2){3}}
\put(4,1){\line(-2,3){4}}

\put(1,7){\line(2,-3){4}}
\put(2,7){\line(1,-2){3}}
\put(3,7){\line(1,-3){2}}

\put(4,7){\line(1,-6){1}}
\put(5,7){\line(0,-6){6}}
\put(6,7){\line(-1,-6){1}}

\put(7,7){\line(-1,-3){2}}
\put(8,7){\line(-1,-2){3}}
\put(9,7){\line(-2,-3){4}}

\put(11,1){\line(2,3){4}}
\put(12,1){\line(1,3){1}}
\put(13,1){\line(0,1){3}}

\put(14,1){\line(1,3){1}}
\put(15,1){\line(0,6){6}}
\put(16,1){\line(-1,3){1}}

\put(17,1){\line(0,1){3}}
\put(18,1){\line(-1,3){1}}
\put(19,1){\line(-2,3){4}}

\put(16,7){\line(2,-3){4}}
\put(17,7){\line(1,-3){1}}
\put(18,7){\line(0,-1){3}}

\put(19,7){\line(1,-3){1}}
\put(20,7){\line(0,-6){6}}
\put(21,7){\line(-1,-3){1}}

\put(22,7){\line(0,-1){3}}
\put(23,7){\line(-1,-3){1}}
\put(24,7){\line(-2,-3){4}}

\put(25,4){\vector(1,0){3}}  
\put(30,1){\line(0,1){6}}
\put(26,3){\makebox(0,0){$\rho$}} 
\put(31,1){\makebox(0,0){$0$}} 
\put(31,4){\makebox(0,0){$1$}} 
\put(31,7){\makebox(0,0){$2$}}

\end{picture}
\caption{ $DL(3^2,3^2)$ includes into $DL(3,3)$.}
\label{dl9933}
\end{figure}
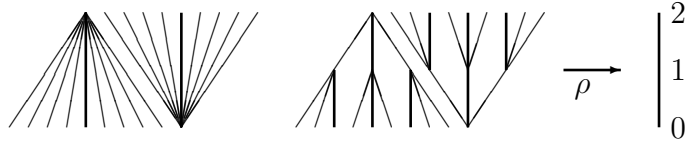
\begin{defn} For any $k$, we can define a map $$up: DL(n,n) \to i(DL(n^k,n^k))$$ by sending
$(p,q)$ to some point $(p',q')$ where 
\begin{itemize}
\item $ht(p')=ht(q')\in k \Z$,  
\item $ht(q')=ht(p') \geq ht(p)=ht(q)$
\item $d_{T_{n+1}}(p,p')=d_{T_{n+1}}(q,q') <  k$.
\end{itemize}
\end{defn}
This map fixes points in $i(DL(n^k, n^k))$ and maps other points ``up'' to the closest level above them in $i(DL(n^k, n^k))$. Of course this $up$ map is not unique but it can always be arranged so that for all $x \in i(DL(n^k, n^k))$ the preimage
has $|up^{-1}(x)|=k$. One could also define this map algebraically but we wish to emphasize the geometry of $DL(n,n)$. Note that the $up$ map is a bounded distance from the identity map on $DL(n,n)$. 

\begin{lemma} Any bijective quasi-isometry $\phi':DL(n^k,n^k) \to DL(n,n)$ is a bounded distance from a map of the form 
$$ i \circ \phi$$
where $\phi: DL(n^k,n^k)  \to DL(n^k,n^k)$ is a $k$ to $1$ quasi-isometry.\end{lemma}
\begin{proof}
Let $\phi= i^{-1} \circ up \circ \phi'$ then $\phi$ is a $k$ to $1$ map and $i \circ \phi = up \circ \phi'$ is a bounded distance from $\phi'$.
Each of $\phi', i$ and $up$ are quasi-isometries so $\phi$ is also a quasi-isometry.
\end{proof}

\section{Boundaries} 
The material in this section can also be found in \cite{EFW} or \cite{FM1} but we reproduce it here for the convenience of the reader.
\subsection{Tree Boundaries}
The oriented tree $T_{n+1}$ has a special class of bi-infinite geodesics, the vertical geodesics. These are lines in $T_{n+1}$ that are coherently oriented. Equivalently they are the lines that project isometrically onto $\R$ under the height map. Similarly, we can define vertical geodesic rays as sub-rays of vertical geodesics. We define the tree boundary to be equivalence classes of vertical geodesic rays. This coincides with the usual notion of tree boundary but our definition allows us to partition the boundary as
$$ \partial_l T_{n+1} \cup \{ \infty \}$$ 
where $\infty$ corresponds to the one class of positively oriented vertical rays. (Recall that each vertex in $T_{n+1}$ has $n$ incoming and one outgoing edges.) We will simply write $\partial T_{n+1}$ to mean $\partial_l T_{n+1}$.

We can identify $\partial T_{n+1}$ with $\Q_n$ by the following procedure. First choose a coherently oriented line in $T_{n+1}$ and label each edge of this line with $0$. Then for each vertex assign labels $0,\ldots,n-1$ to incoming edges.
Now each coherently oriented line $\ell$ in $T_{n+1}$ defines an element of $\Q_n$ by 
$$ \ell \mapsto \sum a_i n^{-i}$$
where $a_i$ is the label of the edge going from height $i$ to $i+1$. Since each line eventually has all edges labeled $0$ this is indeed a map to $\Q_n$. Likewise each element of $\Q_n$ defines a line in $T_{n+1}$. The usual metric on $\Q_n$ has the following geometric interpretation. Consider two points $\xi, \xi'$ on the boundary and two vertical geodesics $\gamma, \gamma'$ emanating from $\xi$ and $\xi'$. These two geodesics meet for the first time at some height $I$. This is precisely the index below which the two series for $\xi, \xi' \in \Q_n$ coincide. Hence $d_{\Q_n}(\xi, \xi')=n^{I}$.

\subsubsection{Clones and measure}\label{candm}
\begin{defn}[clone] Given a point $p \in T_{n+1}$ we define the shadow of $p$ in $\partial T_{n+1}$ to be all $\xi \in \partial T_{n+1}$ that can be reached by vertical geodesics passing through $p$. Note that $diam(C)=n^{h(p)}$. Any subset $C$ of $\Q_n$ that can be defined in this way will be called a clone.  
\end{defn}
We also have a natural measure on $\Q_n$ (Hausdorff measure) but since we will only be concerned with finite unions of clones we can simply define for any clone $C$
$$ \mu(C):= diam(C)$$ and for any finite union of disjoint clones $C= \sqcup A_i$ $$\mu(C):= \sum \mu(A_i).$$ 

\subsection{$DL(n,n)$ boundaries}
As in the tree case we have an orientation on $DL(n,n)$ and a notion of vertical geodesic rays. 
We define two boundaries for $DL(n,n)$. The lower boundary $\partial_l DL(n,n)$ is defined to be equivalence classes of downward oriented vertical geodesic rays whereas the upper boundary $\partial_u DL(n,n)$ is the equivalence classes of upward oriented vertical geodesic rays. As pointed out in \cite{EFW} $\partial_l DL(n,n) \simeq \partial T_{n+1} \simeq \partial_u DL(n,n) \simeq \Q_n$.

Given a point $p\in DL(n,n)$ we can now define two different shadows, the lower shadow $C_l$ and the upper shadow $C_u$. Both of these shadows are clones in $\Q_n$ with the property that $diam(C_u) = 1/diam(C_l)$ or, in terms of measure, $\mu(C_l)\mu(C_u)=1$.

\subsection{Height respecting quasi-isometries}
A height respecting quasi-isometry of $T_{n+1}$ or $DL(n,n)$ is a quasi-isometry that permutes level sets of the height function up to bounded distance such that the induced map on height is simply translation. Theorem \ref{efw} says that all quasi-isometries of $DL(n,n)$ are height respecting.
It is as straightforward computation to see that height respecting  quasi-isometries induce bilipschitz maps of $\partial_l$ and $\partial_u$ (see \cite{EFW}). Likewise, height respecting isometries induce \emph{similarities}, maps that scale distance by a fixed factor, on each of the boundaries. 

\section{Structure of bilipschitz maps of $\Q_n$}

\subsection{Preliminaries}
In this section we use a theorem of Cooper to analyze the structure of surjective bilipschitz maps of $\Q_n$. The definitions and theorems used in this section can be found in \cite{FM1}.

\begin{defn} A map $\phi : \Q_n \to \Q_n$  is said to be \emph{measure linear} on a clone $C \subset \Q_n$ if there exists some $\lambda$ such that for all $A \subset C$ 
$$\frac{\mu(\phi(A))} {\mu(A)} = \lambda.$$
\end{defn}
The following theorem applies to general clone cantor sets $C,C'$ but in our case $C,C'$ will be clones in $\Q_n$ as defined in Section \ref{candm}.\\

\noindent{}{\bf Theorem 10.6 (Cooper)}\cite{FM1} Suppose that $C,C'$ are clones Cantor sets and that $\phi$ is a bilipschitz map of $C$ onto a clopen subset of $C'$. Then there is a clopen $A$ in $C$ such that the restriction $\phi\vert_A$ of $\phi$ to $A$ is measure linear.\\

As pointed out by Cooper in \cite{FM1}, a clopen in $C\subset \Q_n$ is a finite union of clones.  To prove Proposition \ref{measlinclone} below we simply need to show that for a surjective bilipschitz map $\phi: \Q_n \to \Q_n$ the image of a clone is a finite union of clones.  Then we can apply Cooper's theorem.



\begin{prop}\label{measlinclone}
Given a bilipschitz map $\phi: \Q_n \to \Q_n$ there is a clone $A \subset \Q_n$ so that $\phi \vert_A$ is measure linear.
\end{prop}
\begin{proof}
We show that for any clone $C$ the image $\phi (C)$ is a finite union of clones.\\ 

\noindent{}{\bf Claim} $B\subset \Q_n$ is a finite union of clones if and only if  $$sep(B):=\inf \{ d_{\Q_n}(x,y) \mid x \in B , y \notin B \}  > \e$$ for some $\e > 0$. \\

\emph{Proof of Claim.} If $B$ is a finite union of clones then $sep(B)> diam(D_{min})$ where $D_{min}$ is the smallest clone in the union. 
Conversely if $sep(B)> \e$ then consider clones  $D$ with size $diam(D) <  \e$. If $D\cap B \neq \emptyset$ then $D \subset B$ otherwise there would be $x,y \in D$ with $x\in B, y \notin B$ with $d(x,y) < diam(D) <  \e$ contradicting that 
$sep(B)> \e$.\\

 Now given $y \notin \phi(C)$ there is some $z \notin C$ such that $\phi(z)=y$. Since $C$ is a clone there exists $\e_C$ such that $d(z,C)> \e_C$ for all $z \notin C$. 
Then since $\phi$ is bilipschitz  $d(y,\phi(C)) \geq \frac{1}{K} \e_C.$ In particular $sep(\phi(C))>1/K \e_C$. This says that $\phi\vert_C$ maps onto a clopen of $C' \subset \Q_n$.
Now we can apply Cooper's Theorem.
\end{proof}

\subsection{Zooming}\label{tukiasection}
The following proposition was suggested by Souto. 
\begin{prop}\label{Tukiaprop}
If $\phi: \Q_n \to \Q_n$ is a bilipschitz map then there exists a sequence $\phi_i$ of similarities of $\Q_n$ such that 
$$\lim_{i \to \infty} \phi_i^{-1} \phi \phi_i = \bar{\phi}$$
and $\bar{\phi}$ is measure linear on all of $\Q_n$.
\end{prop}
\begin{proof}
Since $\phi$ is bilipschitz then by Proposition \ref{measlinclone} there is a clone $A \subset \Q_n$ such that $\phi \vert_A$ is measure linear. 
We identify $\Q_n\simeq \partial T_{n+1}$, bilipschitz maps  of  $\Q_n$ with height respecting quasi-isometries of $T_{n+1}$ and  similarities of $\Q_n$ with height respecting isometries. In particular if $\phi$ is a quasi-isometry of $T_{n+1}$ then we also write $\phi$ for the induced  boundary map. 

Let $p \in T_{n+1}$ and $\phi_i$ a sequence of isometries of $T_{n+1}$ such that  $$\phi_i(p) \to x \in A$$ and such that there is some vertical geodesic in $\ell \subset T_{n+1}$ and $M\geq 0$ such that $d(\phi_i(p),\ell)\leq M$. 
We can pick such $\phi_i$ since the group of height-respecting isometries acts cocompactly on $T_{n+1}$. Then the set $\{(\phi_i^{-1} \phi \phi_i)(p)\}$ is bounded in $T_{n+1}$ and so the quasi-isometries $\phi_i^{-1} \phi \phi_i$ converge to some quasi-isometry $\bar{\phi}$. This convergence can also be seen on the boundary.  

Now for any clone $C \subset \Q_n$ there is some $j_C$ such that if $i>j_C$ then $C \subset \phi_i^{-1}(A)$ so that $\phi_i^{-1}\phi \phi_i$ is measure linear on $C$ if $i>j_C$. By Lemma \ref{juanlemma2} below, there is some $i_C$ such that for $i>i_C$
$$\phi_i(C)=\phi(C).$$ 
Therefore $\phi$ is measure linear on all of $\Q_n$. \end{proof}
%
%
\begin{lemma}\label{juanlemma2}
If $\phi_i: \Q_n \to \Q_n$ is a sequence of bilipschitz maps that converge to $\phi$ then for all 
clones $C \subset \Q_n$ there exists $i_C$ such that $\phi_i(C)=\phi(C)$ for $i \geq i_C$.
\end{lemma}
\begin{proof}
Since $ \phi_i \to \phi$, there exists $i_\e$ such that for all $x \in C$ 
we have $d(\phi(x),\phi_i(x)) \leq \e$ 
therefore $\phi_i(x) \in \phi(C)$.
Suppose $y_i \in \phi(C)\setminus \phi_i(C)$ then $y_i \to y$. But 
$$y=\phi(x)=\lim_{i \to \infty} \phi_i(x) $$
Since $d(y_i, \phi_i(x)) \to 0$ then by the same argument we have that $y_i \in \phi(C)$ for $i> i_C$.
\end{proof}

\begin{prop}\label{powersofn}
A map $\bar{\phi}: \Q_n \to \Q_n$ that is onto and measure linear on all of $\Q_n$ has measure linear constant 
$$\lambda= r_1^{j_1}\ldots  r_i^{j_i}$$
where $j_l \in \Z$ and $r_l$ is prime divisor of $n$. In particular if $n=r^i$ where $r$ is prime then $\lambda=r^j$ for some $j \in \Z$. 
\end{prop}
\begin{proof}
Pick any $C \subset \Q_n$. Then $\bar{\phi}(C)$ is a finite union of clones (since $\bar{\phi}$ is bilipschitz). Of course any finite union of clones can also be written as a finite union of disjoint clones. Each clone $B$ has $\mu(B)=diam(B)=n^i$ for some $i$. Therefore $\lambda$ is a finite sum of powers of $n$. 
Now $\bar{\phi}$ has an inverse $\bar{\phi}^{-1}$ that is also measure linear with measure linear constant $1/\lambda$.
So both $\lambda$ and $1/\lambda$ must be finite sums of powers $n$.  The only possibilities are $\lambda= r_1^{j_1}\ldots  r_i^{j_i}$ since $\lambda$, when written as a fraction, must have denominator $n^l$ for some $l \in \mathbb{N}$. 
\end{proof}

\begin{cor}\label{keycor} If $\phi:DL(n,n) \to DL(n,n)$ is a quasi-isometry then there exists a quasi-isometry $\bar{\phi}$ such that $\bar{\phi}_l$ and $\bar{\phi}_u$ are both measure linear on all of $\Q_n$.
\end{cor}
\begin{proof} We modify the proof of Theorem \ref{Tukiaprop} and repeat the zooming argument twice. 
The map $\phi$ has lower boundary map $\phi_l$ that is bilipschitz. Let $A_l$ be a clone on which $\phi_l$  is measure linear as provided by Proposition \ref{measlinclone}. Pick $\phi_i$ to be isometries of $DL(n,n)$ such that $\phi_i(p) \to x \in A_l$ and such that $\phi_i(p)$ stays a bounded distance from some vertical geodesic in $DL(n,n)$. Then by similar arguments as in Theorem \ref{Tukiaprop} $\phi_i^{-1} \phi \phi_i \to \hat{\phi}$ where
$\hat{\phi}_l$ has boundary maps $\hat{\phi}_l$ and $\hat{\phi}_u$ with $\hat{\phi}_l$ measure linear on all of $\Q_n$. 
Now $\hat{\phi}$ has upper boundary map $\hat{\phi}_u$ which is bilipschitz and hence is measure linear on some clone $\hat{A}_u$.  Repeating the above procedure with a sequence of isometries $\hat{\phi}_i$ that now zooms into a point $y \in  \hat{A}_u$ we get a new map $\bar{\phi}$ that has boundary maps $\bar{\phi}_l$ and $\bar{\phi}_u$. The upper boundary map $\bar{\phi}_u$ is measure linear by construction and the lower boundary map $\bar{\phi}_l$ is also measure linear since the lower boundary maps of $\hat{\phi}_i$ were all similarities.
\end{proof}

\begin{prop}
If $\phi$ is an $m$ to $1$ quasi-isometry then $\bar{\phi}$ as constructed in Corollary \ref{keycor} is bounded distance from an $m$ to $1$ quasi-isometry. 
\end{prop}
\begin{proof}
Recall that $$\bar{\phi} = \lim_{i \to \infty } \phi_i^{-1} \phi  \phi_i$$ where $\phi_i, \phi_i^{-1}$ are isometries and hence  $1$ to $1$ maps. A composition of an $m$ to $1$ map with a $1$ to $1$ map is again an $m$ to $1$ map and a limit of $m$ to $1$ maps is also an $m$ to $1$ map. 
\end{proof}

\section{Uniformly finite homology}\label{UFH}

Uniformly finite homology was first introduced in \cite{S} but we do not need the explicit construction here so we refer the reader to \cite{S,BW} for more details. Here we will only introduce the zeroth uniformly finite homology group.
 
 \begin{defn} A metric space $X$ is a \emph{uniformly discrete bounded geometry} (UDBG) space if there exists $\e >0$ such that $d(x,y) > \e$ for all $x \neq y$  and for all $r>0$ there is a bound $M_r$ on the size of any $r$ ball. 
 \end{defn}
For a UDBG space $X$, let $C_0^{uf}(X)$ denote the vector space of infinite formal sums of the form 
$$c= \sum_{x \in X} a_x x  \quad (a_x \in \Z)$$
for which there exists $M_c>0$ such that $|a_x|<M_c$ for all $x \in X$. Let $C_1^{uf}(X)$ denote the vector space of infinite formal sums of the form 
$$ c= \sum_{x,y \in X} a_{(x,y)}(x,y)$$
for which there exists $M_c,R_c >0$ such that $|a_{(x,y)}|< M_c$ and $a_{(x,y)}=0$ if $d(x,y)> R_c$. 
The boundary map is defined by setting
\bea
\partial: C_1^{uf}(X) &\to& C_0^{uf}(X)\\
(x,y) & \mapsto& y-x
\eea
and extending by linearity. The zeroth homology group is then 
$$ H_0^{uf}(X) = C^{uf}_0(X)/\partial(C_1^{uf}(X)).$$
All Cayley graphs are UDBG spaces.
Futhermore, this homology group is a quasi-isometry invariant so if $X$ and $Y$ are quasi-isometric then $H^{uf}_0(X)\simeq H^{uf}_0(Y)$. In particular, for a given group, the uniformly finite homology group does not depend on generating set. 
\begin{defn} Any subset $S\subset X$ defines a class $[S] \in H_0^{uf}(X)$, which is the class of the chain $\sum_{x\in S}x$. We call $[X]$ the \emph{fundamental class} of $X$ in $H_0^{uf}(X)$.
\end{defn}
Using uniformly finite homology, Whyte developed in \cite{W} a test to determine 
when a quasi-isometry between UDBG spaces is a bounded distance from a
bijection.

\begin{thm}\label{Whyte}\cite{W} Let $f:X \to Y$ be a quasi-isometry between
UDBG-spaces. Then there exists a bijective map 
a bounded distance from $f$ if and only if $f_{\ast}([X])=[Y]$ where $$f_{\ast}([X])=\left[\sum_{y\in f(X)}|f^{-1}(y)| y\right].$$

\end{thm}
Another fact that we will use is that two maps $f,g$ that are a bounded distance apart have $$f_*{([X])}=g_*{([X])}.$$
The following theorem allows us to check when a chain $c$ represents the class of the zero chain in $H^{uf}_0(X)$. 
We use this theorem in Proposition \ref{notkto1}.
\begin{thm}\label{whytethm2}\cite{W}
Let $X$ be a UDBG-space, and $c=\sum_{x\in X}a_xx \in C_0^{uf}(X)$. Then we have $[c]=0 \in H_0^{uf}(X)$ if and only if there exist an $r$ such that for any F\o lner sequence $\{S_i\}$, 
\bea
        \left| \sum_{x\in S_i}a_x \right|=O(|\partial_r S_i|).
\eea 
\end{thm}
{\bf Remark.} If $f: X \to X$ is $m\ to\ 1$ then $f_*([X])=m[X]$. If $f$ is a bounded distance from an $m$ to $1$ map then we also have $f_*([X])=m[X]$. In case $m$ is not an integer we will use this condition as a the definition.

\section{From boundary maps to interior maps}\label{coords}\label{constructQI}
In this section we construct an explicit quasi-isometry $\psi$ of $DL(n,n)$ from two boundary maps $\bar{\phi}_l$ and $\bar{\phi}_u$.
First we need to relate $DL(n,n)$ to $\Q_n \times \Q_n \times \R$. 
For this purpose we define two maps: 
$$\pi: \Q_n \times \Q_n \times \R  \to DL(n,n),$$
$$\bar{\pi}: DL(n,n) \to \Q_n \times \Q_n \times \R. $$
The map $\pi$ collapses sets of diameter $n^{\lfloor h\rfloor}$ at height $\ h $ to a single point at height $\lfloor h \rfloor$.
This is a map onto the vertices of $DL(n,n)$. (We could have also defined this map so that it also mapped onto edges but for our purposes this is sufficient.)
The map $\bar{\pi}$ is chosen so that $\pi \bar{\pi}= id$. 
\begin{figure}[h]
\center
\setlength{\unitlength}{.5in} 
\begin{picture}(9,4.6)(1,.5) 
\linethickness{1pt} 

\put(1,1){\line(1,2){1.78}}
\put(2,1){\line(0,1){2}}
\put(3,1){\line(-1,2){1}}

\put(1,1){\circle*{.15}}
\put(2,1){\circle*{.15}}
\put(3,1){\circle*{.15}}

\put(2,3){\circle*{.15}}
\put(2.78,4.6){\circle*{.15}}

\put(4,3){\vector(2,0){2}}  
\put(5,3.25){\makebox(0,0){$\bar{\pi}$}} 

\multiput(7,1)(0,.35){11}{\line(0,1){0.1}} 
\multiput(7.3,1)(0,.35){11}{\line(0,1){0.1}} 
\multiput(7.6,1)(0,.35){11}{\line(0,1){0.1}} 
\multiput(8,1)(0,.35){11}{\line(0,1){0.1}} 
\multiput(8.3,1)(0,.35){11}{\line(0,1){0.1}} 
\multiput(8.6,1)(0,.35){11}{\line(0,1){0.1}} 
\multiput(9,1)(0,.35){11}{\line(0,1){0.1}} 
\multiput(9.3,1)(0,.35){11}{\line(0,1){0.1}} 
\multiput(9.6,1)(0,.35){11}{\line(0,1){0.1}} 

\put(7.3,1){\circle*{.15}}
\put(8,1){\circle*{.15}}
\put(9,1){\circle*{.15}}
\put(7,3){\circle*{.15}}
\put(8.6,4.6){\circle*{.15}}







%

\end{picture}
\caption{A possible $\bar{\pi}$ shown only in projection onto $T_3$ and $\Q_3 \times \R$.}
\end{figure}
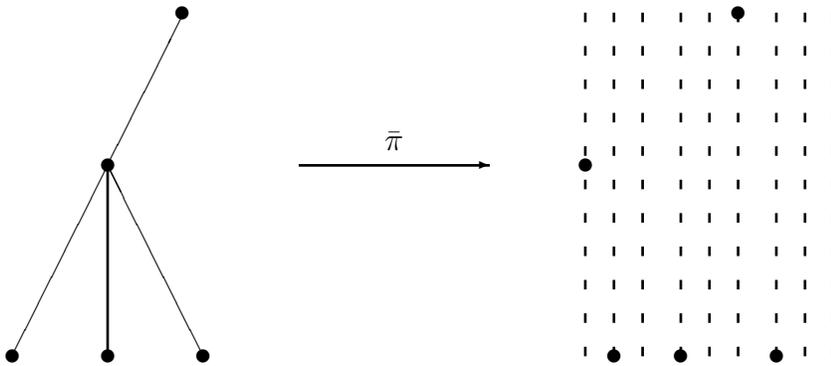
\noindent{}We can endow $\Q_n \times \Q_n \times \R$ with a `sol' like metric so $\pi, \bar{\pi}$ are quasi-isometries.
Specifically, for two points that differ only in the $\R$ coordinate we define their distance to be that difference.
For two points at the same height $(x,y,t), (x',y',t) \in \Q_n \times \Q_n \times \R$ we define $$d((x,y,t),(x',y',t))=n^{-t}d_{\Q_n}(x,x') + n^{t}d_{\Q_n}(y,y').$$
For points $p,q$ at different heights we consider all finite sequences $\{p_i\}$ connecting $p$ and $q$ where subsequent terms are either at the same height or differ only in the $\R$ coordinate and define the distance to be 
 $$d(p,q)= \inf  \sum d(p_i,p_{i+1}).$$
\begin{defn}
Given two bilipschitz maps $\bar{\phi}_l, \bar{\phi}_u$ of $\Q_n$, we define a quasi-isometry $$\psi: DL(n,n) \to DL(n,n)$$
by 
$$\psi:={\pi}\circ (\bar{\phi}_l \times \bar{\phi}_u  \times Id) \circ \bar{\pi}.$$
\end{defn}
%
%


\begin{lemma} Pick two clones $C_l, C_u \subseteq \Q_n$  with $\mu(C_l) \mu(C_u) >1$. Then the image 
$$S_{C_l,C_u}=\pi ( C_l \times C_u \times [ -\log_n (\mu(C_u)), \log_n(\mu(C_l)) ])$$
defines a box of height $H:=|\log_n (\mu(C_u)\mu(C_l))|$ with $\mu(C_l)\mu(C_u)$ vertices at each height. 
\end{lemma}
\begin{proof}
This follows from the definition of $\pi$.
Recall that $|S_{C_l,C_u}|= (H + 1)n^H$ while $|\partial S_{C_l,C_u}|= 2n^{H+1}$. So if we chose a sequence $S_i$ of such boxes with height $H_i$ increasing then $\{S_i\}$ is a F\o lner sequence. 
\end{proof}


\begin{lemma}\label{errest} Let $S_{C_l,C_u}$ be a box as above with height $H\gg \log_n K$. Then, for $r=\log_n K$,
$$ \frac{1}{\lambda_l \lambda_u} | S_{C_l,C_u}|  -  |\partial_r S_{C_l,C_u}| \leq  \sum_{x \in S_{C_l,C_u}} |\psi^{-1}(x)| \leq \frac{1}{ \lambda_l \lambda_u} | S_{C_l,C_u}| + K^2 |\partial_r S_{C_l,C_u}|.$$
\end{lemma}
\begin{proof}
First we claim that since $\bar{\phi}_l$ and $\bar{\phi}_u$ are $K$ bilipschitz then $$\bar{\phi}_l^{-1}(C_l) = \sqcup A_i\ \textrm{ and }\ \bar{\phi}_u^{-1}(C_u) = \sqcup B_i$$ where $A_i$ and $B_i$ are clones of size $\mu(A_i) \geq \frac{1}{K} \mu(C_l)$ and $\mu(B_i) \geq \frac{1}{K} \mu(C_u)$. We already know that the preimages of $\bar{\phi}_l$ and $\bar{\phi}_u$ are clones. To get an estimate on the size note that 
$$n \mu(A_i) = sep(A_i)=sep(C_l) \geq \frac{1}{K} n \mu(C_l)$$
and similarily for the $B_i$. Without loss of generality we can assume that all $A_i$ have the same size $m_A$ and all the $B_i$ have size $m_B$. In particular there are $\frac{\mu(C_l)}{ \lambda_l  m_A}$ many clones in $\bar{\phi}_l^{-1}(C_l)$ and 
$\frac{\mu(C_u)}{ \lambda_l  m_B}$ many clones in $\bar{\phi}_u^{-1}(C_u)$. This implies that, for a fixed height $t$, there are 

$$\frac{\mu(C_l)}{ \lambda_l  m_A} \cdot \frac{\mu(C_u)}{ \lambda_l  m_B}= \frac{\mu(C_l)\mu(C_u)}{ \lambda_l \lambda_u } \cdot  \frac{1}{m_A m_B}$$
many
sets of the form $A_i \times B_j \times \{t\}$ in the preimage of  $\bar{\phi}_l \times \bar{\phi}_u \times id$.
If $ \log_n m_A \geq t \geq - \log_n m_B$  and $t \in \Z$ then $\pi( A_i \times B_j \times \{ t\})$ contains 
$m_A m_B$ 
many vertices but more importantly $$\bar{\pi}(DL(n,n)) \cap A_i \times B_j \times \{t\}$$  contains $m_A m_B$ images of vertices.
For $t$ outside this range this may not be the case and depends on the choice of $\bar{\pi}$.

Now  $\log_n m_A \leq  \log_n(\mu(C_l)) - \log_n K$ and  $- \log_n m_B \geq -\log_n (\mu(C_u))+ \log_n K$ so that 
for height $$t \in  [ -\log_n (\mu(C_u))+ \log_nK, \log_n(\mu(C_l)) - \log_nK ])$$
we have 
$$ \frac{\mu(C_l)\mu(C_u)}{ \lambda_l \lambda_u } \cdot  \frac{1}{m_A m_B} \cdot m_A m_B=  \frac{\mu(C_l)\mu(C_u)}{ \lambda_l \lambda_u } $$
many vertices being mapped into $S_{C_l,C_u}$ at height $t$. But there are $\mu(C_l)\mu(C_u)$ many vertices at this level. Therefore, on average, there are $1/\lambda_l \lambda_u$ many vertices being mapped onto each vertex.

Finally, for $t$ in the range $$ \log_n(\mu(C_l)) \geq  t \geq \log_n(\mu(C_l)) - \log_n{K} \textrm{ or }  \log_n(\mu(C_u)) \leq  t \leq \log_n(\mu(C_u ))  + \log_n{K}$$
 it is possible, depending on the choice of $\bar{\pi}$, that no vertices are mapped to these heights or that as many as $$diam(\phi_l^{-1}(C_l))\cdot diam(\phi_u^{-1}(C_u))$$ vertices are mapped to each height. But $diam(\phi_l^{-1}(C_l)) \leq K \mu(C_l)$ and $diam(\phi_u^{-1}(C_u)) \leq K \mu(C_u)$ so there are at most $K^2 \mu(C_l)\mu(C_u)$ vertices being mapped to each level or in other words on average at most $K^2$ to each vertex. Note that these levels are contained in $\partial_{\log_n K} S_i$.
\end{proof}

\begin{prop}\label{notkto1} Let $X=DL(n,n)$. If $\phi_l, \phi_u$ are measure linear with constants $\lambda_l, \lambda_u$  then $$\psi_{*}(X)\neq k[X] $$ if $k \neq 1/\lambda_l \lambda_u$. 
\end{prop}
\begin{proof}
We will use Whyte's Theorem \ref{whytethm2}. Let $S_i$ be a sequence of boxes in $DL(n,n)$ of height $i$. 
Let $c$ be the chain defined by 
$$c = \sum_{x\in X} a_x x$$
where $a_x=|\psi^{-1}(x)| -k$. 
By Lemma \ref{errest} we have that for $r=\partial_{\log_nK}$ 
$$  \left(\frac{1}{\lambda_l\lambda_u} - k\right)|S_i| - err |\partial_r S_i| \leq  \left| \sum_{x \in S_i} |\psi^{-1}(x)|- k \right|   \leq  \left(\frac{1}{\lambda_l\lambda_u} - k\right)|S_i| +  err |\partial_r S_i|$$
so that unless $1/\lambda_l\lambda_u = k$  we have  $|\sum_{x\in S_i} a_x | \notin O(\partial_r S_i)$ since $|S_i| \notin O(\partial_r S_i)$.
\end{proof}

To conclude the proof of Theorem \ref{thebigone}, we note that since $\psi$ and $\bar{\phi}$ have the same boundary maps they are a bounded distance apart. 
In particular,
$$\psi_*{([X])}= \bar{\phi}_{*}([X]).$$
Now if $\bar{\phi}$ is $k$ to $1$ then $\bar{\phi}_{*}([X])=k[X]$ but, as shown in Proposition \ref{notkto1}, this is impossible.


\begin{thebibliography}{99}




\bibitem[BW1]{S}
 Block, J.,Weinberger, S.,
\textit{Aperiodic tilings, positive scalar curvature and amenability of spaces},
  J. Amer. Math. Soc.  5  (1992),  no. 4, p. 907-918. 

\bibitem[BW2]{BW}
Block, J., Weinberger, S.,
\textit{Large scale homology theories and geometry},
AMPS/IP Studies in Advanced Mathematics, 2 (1997), p. 522-569.


\bibitem[BK]{Kl} 
Burago, D., and Kleiner, B.,
\textit{Separated nets in Euclidean space and Jacobians of bi-Lipschitz maps},
Geom. Funct. Anal. 8 (1998), no. 2, 273-282. 






\bibitem[D]{D}Dymarz, T.,
\textit{Bijective quasi-isometries of amenable groups}, Geometric methods in group theory, 
Contemp. Math., 372, Amer. Math. Soc., Providence, RI, (2005), 181-188. 

\bibitem[EFW1]{EFW} Eskin, A., Fisher, D., and Whyte, K.,
  \textit{Quasi-isometries and rigidity of solvable groups}, 
 Pure Appl. Math. Q.  3  (2007),  no. 4, part 1, 927-947. 
 
\bibitem[EFW2]{EFW2} Eskin, A., Fisher, D., and Whyte, K.,
  \textit{Coarse differentiation of quasi-isometries I: Rigidity for Sol and Lamplighter groups}, to
  appear.

\bibitem[EFW3]{EFW1} Eskin, A., Fisher, D., and Whyte, K.,
  \textit{Coarse differentiation of quasi-isometries II: Rigidity for Sol and Lamplighter groups}, to
  appear.

\bibitem[FM1]{FM1} Farb, B., and Mosher, L., \textit{A rigidity
    theorem for the solvable Baumslag-Solitar groups}, With an
  appendix by Daryl Cooper. Invent. Math. 131 (1998), 419-451.
  
\bibitem[H]{H} de la Harpe, P., \textit{Topics in Geometric Group Theory}, Chicago, 2000







    




















\bibitem[TW]{TW} Taback, J., and Wong, P.,
\textit{The geometry of twisted conjugacy  classes in wreath products},
Preprint, (2008)

\bibitem[Wh]{W} Whyte, K.,
\textit{Amenability, bi-Lipschitz equivalence, and the von Neumann conjecture},
Duke Math. J. 99 (1999), no. 1, 93-112. 

\bibitem[Wo]{Wo} Wortman, Kevin., \textit{A finitely presented solvable
    group with small quasi-isometry group,} Michigan Math. J. 55 (2007), 3-24.


\end{thebibliography}
\end{document}